
\magnification=\magstep1
\baselineskip=16pt plus 1pt minus 1pt

\centerline{\bf Connection between q-deformed anharmonic oscillators and}

\centerline{\bf quasi-exactly soluble potentials}

\bigskip\bigskip\bigskip

\centerline{Dennis Bonatsos $^{1,2}$\footnote*{e-mail:
bonat@ectstar.ect.unitn.it},
C. Daskaloyannis $^3$\footnote{$^+$}{e-mail: daskaloyanni@olymp.ccf.auth.gr}
  and Harry A.
Mavromatis $^4$\footnote{$^{\#}$}{e-mail: facm007@saupm00.bitnet }}

\centerline{$^1$  European Centre for Theoretical Studies in Nuclear
Physics and Related Areas (ECT$^*$)}

\centerline{Strada delle Tabarelle 286, I-38050 Villazzano (Trento), Italy}

\centerline{$^2$ Institute of Nuclear Physics, NCSR ``Demokritos''}

\centerline{GR-15310 Aghia Paraskevi, Attiki, Greece}

\centerline{$^3$ Department of Physics, Aristotle University of Thessaloniki}

\centerline{GR-54006 Thessaloniki, Greece}

\centerline{$^4$ Physics Department, King Fahd University of Petroleum and
Minerals}

\centerline{ Dhahran 31261, Saudi Arabia}

\bigskip\bigskip\bigskip

\centerline{\bf Abstract}

It is proved that quasi-exactly soluble potentials (QESPs) corresponding to
an oscillator with harmonic, quartic and sextic terms,
for which the $n+1$ lowest levels of a given parity can be  determined
exactly, may be approximated by  WKB equivalent potentials corresponding to
deformed anharmonic oscillators of SU$_q$(1,1) symmetry, which have been used
for the description of vibrational spectra of diatomic molecules.
This connection allows for the immediate approximate determination of all
levels of the QESPs.

\vfill\eject

Quantum algebras [1,2] (also called quantum groups), which are deformed
 versions of the usual Lie algebras, to which they reduce when the deformation
parameter is set equal to unity,  have  recently been attracting considerable
attention. The interest in possible physical applications was triggered by the
introduction of the $q$-deformed harmonic oscillator by Biedenharn [3] and
Macfarlane [4] in 1989, although similar oscillators had already been in
existence [5-7]. The common mathematical structure of these various types
of oscillators was shown by several authors (see [8] for a list
of references), using the concept of the generalized deformed oscillator
[9].

One way to clarify the physical content of deformed oscillators is to try
to construct potentials which give, exactly or  approximately, the same
spectrum as these oscillators. Along these lines, WKB equivalent potentials
corresponding to the $q$-deformed oscillator [10], as well as to deformed
anharmonic oscillators [11] found useful in the description of vibrational
spectra of diatomic molecules [12,13] have been constructed, both
numerically [10] and in analytic form [11]. These
potentials give approximately the same spectrum as the deformed oscillators
under discussion.

On the other hand, several quasi-exactly soluble potentials (QESPs) have been
introduced [14-18], for which the exact calculation of the first $n+1$
energy levels is possible, with no information provided for the rest of the
levels. In particular, in the case of the harmonic oscillator with both
quartic and sextic
anharmonicities [14], the exact calculation of the first
$n+1$ energy levels {\sl of a given parity}  is possible.

It is therefore interesting to check if there is any relation between the
WKB equivalent potentials (WKB-EPs) approximating the behavior of the deformed
oscillators on the one hand and the quasi-exactly soluble potentials on the
other. If such a link exists, as we shall demonstrate in this letter,
one can fix the parameters of a WKB potential
equivalent to an appropriate oscillator so that the potential agrees,
up to the sextic term, with a given quasi-exactly soluble potential.
If the approximations involved (WKB method, truncation of the WKB
potential to sixth order only) are good, one can subsequently  use
the levels of the relevant oscillator in order to approximate the
unknown higher levels of the quasi-exactly soluble potential. In the case
of the harmonic oscillator with both quartic and sextic anharmonicities
[14], the levels of the opposite parity can also be approximated in this
way.

The q-numbers are defined as
$$[x]= {q^x -q^{-x}\over q-q^{-1}}.\eqno(1)$$
In the case that $q=e^{\tau}$, where $\tau$ is real, they can be
written as
$$[x]={\sinh (\tau x) \over \sinh (\tau)}, \eqno(2)$$
while in the case that $q$ is a phase ($q=e^{i\tau}$, with $\tau$
real) they take the form
$$[x]={\sin (\tau x)\over \sin (\tau)}. \eqno(3)$$
It is clear that in both cases the q-numbers reduce to the usual
numbers as $q\rightarrow 1$ ($\tau\rightarrow 0$).

The q-deformed harmonic oscillator [3,4] is determined by
 the creation and annihilation operators $a^+$ and $a$, which
satisfy the  relations:
$$ a a^+ -q^{\mp 1} a^+ a = q^{\pm N}, \quad
[N, a^+] = a^+, \quad [N, a]= -a , \eqno(4)$$
where $N$ is the number operator.
The Hamiltonian of the q-deformed oscillator is
$$ H= {\hbar \omega \over 2} (a a^+ + a^+ a), \eqno(5)$$
and its   eigenvalues are [3,4]
$$ E(n) = {\hbar \omega\over 2} ([n]+[n+1]). \eqno(6)$$

{}From Eq. (6) it is clear that the spectrum of the q-deformed
harmonic oscillator is not equidistant. For real $q$ the
spectrum increases more rapidly than in the classical case,
while for $q$ complex it increases less rapidly than
in the classical case, i.e. it gets squeezed.

Using two q-deformed oscillators described by the operators
$a_1,a^+_1,N_1$, and $a_2,a^+_2,N_2$, the generators of the
quantum algebra SU$_q$(1,1) are written as [19,20]
$$K_+=a_1^+a_2^+, \quad K_=a_1a_2, \quad
K_0={1\over 2}\left( N_1+N_2+1\right), \eqno(7)$$
and satisfy the commutation relations
$$ \left[ K_0,K_\pm \right] = \pm K_\pm, \quad
\left[ K_+,K_- \right] = - \left[ 2K_0 \right],\eqno(8)$$
where the square bracket in the rhs of the last equation is a q-number
as defined in eq. (1). In ref. [12]
this symmetry has been used for the description of
vibrational molecular spectra. The q-deformed anharmonic oscillator
used there is a q-generalization of the anharmonic oscillators used
in the usual Lie algebraic approach to molecular spectroscopy.
(See ref. [12] for detailed references.)

The energy  spectrum of the anharmonic oscillator with SU$_q$(1,1)
symmetry, in the case of complex $q$, is
given by [12]
$$ E_n = E_0' -A {{\sin(\tau(n-N/2))\sin(\tau(n+1-N/2))} \over
{\sin^2 (\tau)}} , \eqno(9)$$
where
$$N=2n_{max} \quad\quad \hbox{or} \quad\quad N=2n_{max}+1,\eqno(10)$$
with $n_{max}$ corresponding to the last level before the
dissociation limit.
For $\tau \to 0$ a Morse [21] or modified P\"oschl-Teller [22,23] spectrum
is obtained.
The WKB equivalent potential in this case is [11]
$$ V(x)= V_{min} +
{A\over 4} \left({ \tau \sin (N\tau )\over  \sin^2 \tau} \right)^2
u^2$$
$$ \Big[1 - {2\over 3} {\tau^2 \cos (N\tau )\over \sin^2 \tau}u^2
+{1\over 45} { \tau^4 (23\cos^2(N\tau)-6)\over
\sin^4\tau}u^4 $$
$$- {2\over 315}
{ \tau^6 ( 67\cos^2(N\tau)-36)\cos(N\tau) \over \sin^6\tau}u^6+\cdots
\Big], \eqno(11)$$
where
$$ V_{min} = E_0' - A {\cos (\tau) -\cos (\tau N)\over 2 \sin ^2 \tau},
\eqno(12)
$$
and
$$u={ \sqrt{2mA} \over \hbar} x. \eqno(13)$$
It can be easily seen that for $\tau \to 0$ Eq. (11) reduces to the
Taylor expansion of the modified P\"oschl-Teller potential
$$ V(x) = V_{min} + {A\over 4} u^2 \left[ 1-{2\over 3} u^2 +{17\over 45} u^4
-{62\over 315} u^6 + \ldots \right], \eqno(14)$$
where $u$ is given by eq. (13),
while  the modified P\" oschl--Teller potential [22,23] in closed form is
$$ V(x)=V_{min} + {A N^2 \over 4}
\tanh^2\left({ \sqrt{2mA}x \over \hbar}\right).\eqno(15)$$
 This means that the WKB-EP of a system with
$SU_q(1,1)$ symmetry is a deformation of the modified P\"oschl--Teller
potential.

On the other hand, it is known [14,15] that the potentials
$$ V(x) = 8 a^2 x^6 + 8 a b x^4 + 2 [b^2 -(2k+3) a ] x^2, \eqno(16)  $$
with $k=2n+r$, where $n=0$, 1, 2, \dots and $r=0$, 1 are quasi-exactly
soluble. The meaning of this term is the following: for these potentials
one can construct exactly the first $n+1$ levels with parity $(-1)^r$.
(The extra factors of 2 in eq. (16) in comparison to refs [14,15] are due
to the fact that we use the usual form of the Schr\"odinger
operator
$$ H= - {\hbar^2 \over 2 m} {d^2 \over d x^2 } + V(x),\eqno(17)$$
while in refs [14,15] the form
$$ H = -{\hbar^2 \over m} {d^2 \over d x^2} + V(x) \eqno(18)$$
is used.)

We wish to check if the WKB-EP for the deformed oscillator mentioned
above corresponds, up to terms of $x^6$, to quasi-exactly soluble potentials
of the form of eq. (16). For reasons of convenience in eq. (13) we set
$\hbar = m=1$.
Comparing the coefficients of $x^2$, $x^4$ and $x^6$ in eqs (11) and (16)
we get respectively the following equations
$$ 2 ( b^2 -(2k+3) a) = {A^2\over 2} {\tau^2 \sin^2(N\tau) \over
\sin^4\tau}, \eqno(19)$$
$$ 8 a b = -{2\over 3} A^3 {\tau^4 \sin^2(N\tau) \cos(N\tau) \over
\sin^6\tau}, \eqno(20)$$
$$ 8 a^2 = {2\over 45} A^4 {\tau^6 \sin^2(N\tau) (23\cos^2(N\tau)-6)\over
\sin^8 \tau}.\eqno(21) $$
We can consider $a^2=1$ without loss of generality.
For the ground state to be normalizable one should have $a\geq 0$ [16],
which implies in the present case $a=1$. Having in mind that
in eq. (9) one must have $A>0$ in order to get an increasing spectrum
(this fact can be seen also in the realistic applications to molecular
spectra given in ref. [12]), from eq. (21) one finds
$$ A= \sqrt{6\sqrt{5}} {\sin^2\tau \over \tau^{3/2} \sqrt{\sin(N\tau)}
\root 4 \of {23\cos^2(N\tau)-6} } .\eqno(22) $$
Because of the symmetry $q\leftrightarrow q^{-1}$ characterizing the
q-deformed oscillators, it suffices to consider $\tau>0$. The following
conditions should then be satisfied
$$ \sin(N\tau) >0, \eqno(23)$$
$$ \cos^2(N\tau)>6/23.\eqno(24)$$
Then eq. (20) gives
$$ b = -\sqrt{ {15\sqrt{5}\over 2} } {\sqrt{\sin(N\tau)} \cos(N\tau) \over
\sqrt{\tau} (23\cos^2(N\tau)-6)^{3/4} },\eqno(25) $$
while eq. (19) gives
$$ 2k+3 = - {3\sqrt{5} \over 2} {\sin(N\tau) (18\cos^2(N\tau) -6) \over
\tau (23\cos^2(N\tau)-6)^{3/2} }.\eqno(26)$$
The fact that $k=2n+r$ has to be positive gives the additional condition
$$ \cos^2(N\tau) < 1/3.\eqno(27) $$
{}From eq. (25) it is clear that $b$ can be either positive (when
$\cos(N\tau)<0$) or negative (when $\cos(N\tau)>0$).

For reasons of convenience, let us consider the case $b>0$ first.
In this case the conditions $\sin(N\tau)>0$ and $\cos(N\tau)<0$ imply that
we can limit ourselves to the region $\pi > N\tau >\pi/2$. The conditions
(24) and (27) then imply  that $ 2.1863 + 2 \pi l > N\tau > 2.1069 + 2 \pi l$,
with $l=0$, 1, 2, \dots. Thus though we can find an infinity of values for
$N$, as we shall see it suffices to consider $l=0$.

We first wish to check the accuracy of our approach, which already involves
two major approximations: the WKB approximation and the omission in the
WKB-EP of terms higher than $x^6$. In order to achieve this, we will compare
the results given by this method to the exact results given in
ref. [14] for the case $n=1$, $r=0$, in which $2k+3=7$. From eq. (26) one
then sees that a possible solution is $N=151$, $\tau=0.0144503$, which
gives $A=0.4343473$, $b=12.589097$, while from eq. (12) (for $V_{min}=0$) one
 has $E'_0= 1636.8943$.
We see therefore that eqs (19)--(21) are indeed satisfied and the potentials
of eqs (11) and (16) (up to the $x^6$ term) are identical, being
$$ V(x) = 303.02 ( x^2 + 0.3324 x^4 + 0.02640 x^6).\eqno(28) $$
The two solutions which can be obtained exactly are [14]
$$E_0 = 3 b -2\sqrt{b^2+2} = 12.4307, \qquad E_2 = 3b+ 2\sqrt{b^2+2} =
63.1039, \qquad E_2-E_0 = 50.6731,\eqno(29) $$
while eq. (9) gives the complete spectrum, including
$$E_0 = 12.3805, \qquad E_2 = 63.0584, \qquad E_2-E_0= 50.6779. \eqno(30) $$
We remark that the agreement between the exact results of ref. [14] and the
predictions of eq. (9) is excellent, implying that in this case:

i) the WKB method is accurate,

ii) the omission of terms higher than $x^6$ is a good approximation.

An additional check for the case $n=3$, $r=0$ is also made. In this case
$2k+3=15$. One possible solution of eqs (19)-(21) in this case is given by
$N=325$, $\tau=0.00671384$, $A=0.2960795$, $b= 18.469158$, $E'_0=5168.941$.
The potential is
$$V(x) = 652.20 (x^2 + 0.2265 x^4 + 0.01227 x^6).\eqno(31) $$
Eq. (9) gives the levels
$$ E_0= 18.1126, \qquad E_2 = 91.3482, \qquad E_4=165.8771, \qquad
E_6 = 241.6456,\eqno(32)$$
while, following the procedure of ref. [14], one sees that in this case
the exact energy levels are the roots of the equation
$$ E^4 -28 b E^3 + (254 b^2 -240) E^2 + (-812 b^3 + 2592 b) E +
(585 b^4 -4656 b^2 +2880) =0,\eqno(33) $$
given by
$$ E_0=18.1429, \qquad E_2=91.3783, \qquad E_4=165.913, \qquad E_6=241.703.
  \eqno(34)$$
We remark that excellent agreement between the approximate and the exact
results is again obtained.

The agreement remains equally good if more levels are considered.
For the  case $n=9$, $r=0$, in which $2k+3=39$, one possible solution
of eqs (19)-(21) is given by $N=399$, $\tau=0.00545864$, $A=0.2703882$,
$b=21.275801$, $E_0'= 7126.0336$. The potential is
$$ V(x)= 827.43 (x^2 + 0.2057 x^4 + 0.009669 x^6).$$
The energy levels are compared to the exact ones in Table I.

We turn now to the case with $b<0$. The conditions $\sin(N\tau)>0$ and
$\cos(N\tau)>0$ imply that we should limit ourselves to the region
$\pi/2 > N\tau >0$. Then the conditions (24) and (27) imply that
$ 1.0347 > N\tau > 0.9553$. (We again ignore terms of $2 \pi l$.)

We again consider the case $n=1$, $r=0$, $2k+3=7$. One possible solution is
given by $N=61$, $\tau=0.0157377$, $A=0.4538508$, $b=-12.108743$, $E'_0=
390.66689$.
The relevant potential is
$$ V(x) = 279.095 (x^2 - 0.3471 x^4 + 0.02866 x^6).\eqno(35) $$
The exact solutions given by ref. [14] are
$$ E_0 = -60.7083, \qquad E_2 = -11.9441,\eqno(36)$$
while eq. (9) gives
$$ E_0 = 11.7456, \qquad E_2= 57.3764.\eqno(37) $$

The reasons for this dramatic failure are well understood. The potential
of eq. (35) has a central well in the middle, having its minimum at
$x=0$, $E=0$, plus two symmetrically located  wells, on either side of the
central one, with their minima at $E<0$. The method of ref. [14]
gives the lowest two energy eigenvalues of the system, which in this case are
the lowest two levels in the side wells. The WKB method used in ref. [11],
however, is known to be valid only for small values of $x$, i.e. in the
area of the central well. Therefore in this case eq. (37) gives the lowest
two levels of the central well, which in this case are {\sl not} the lowest
two levels of the system.

We therefore conclude that a correspondence between the SU$_q$(1,1) WKB-EP
and the quasi-exactly soluble potentials of eq. (16) can be made only
for $b>0$, in which case both methods will give the levels corresponding
to a well around $x=0$. In this case, as we have already seen, the
approximation is very good, providing us with the following method:

 Given a quasi-exactly soluble potential for which only the $n+1$ lowest
levels of parity $(-1)^r$ can be found exactly, we can choose appropriately the
parameters of an SU$_q$(1,1) WKB-EP, so that the two potentials are
identical up to order $x^6$. Then using eq. (9), i.e. the known eigenvalues
of  the WKB-EP, we can approximate the levels of the same parity lying above
the first $n+1$ known levels, as well as the levels of opposite parity,
for which the method of ref. [14] gives no information {\sl for the
potential under discussion}. (It is worth noting that changing the
values of $n$ and $r$ changes the potential.)

Using a similar procedure one can easily see that no connection
is possible between the WKB-EPs [11] of the q-deformed harmonic oscillator
and the QESPs. In the case of the Q-deformed harmonic oscillator [5--7],
the WKB-EP of which has been given in [24], no connection is possible either.
The same is true for the modified P\"oschl-Teller potential, the Taylor
expansion of which is given by eq. (14).

The following comments on the results can now be made:

i) The quasi-exactly soluble potentials are known to be related to the SL(2)
[15] and deformed SU(2) [25] symmetries. Among the several oscillators
mentioned here, the only one for which the connection to QESPs was possible
is the oscillator with SU$_q$(1,1) symmetry. Apparently the H$_q$(4) symmetry
of the q-deformed oscillator is ``not  enough'' to guarantee such a
connection.

ii) The QESPs of eq. (16) are characterized by 3 parameters ($a$, $b$, $k$).
The SU$_q$(1,1) potential is also characterized by 3 parameters ($A$, $\tau$,
$N$), while the WKB-EPs of the q-oscillator and Q-oscillator are characterized
 by 2 parameters ($\omega$, $\tau$ and $\omega$, $Q$ respectively) and the
modified P\"oschl-Teller potential is characterized by only one parameter
($A$). It is therefore not surprising that a connection is possible only for
the case in which the number of parameters of the two potentials to be related
is the same.

In summary, we have proved that the quasi-exactly soluble potentials
corresponding to a harmonic oscillator with quartic and
sextic anharmonicities  can be approximated by the WKB equivalent potentials
corresponding to a deformed anharmonic oscillator with SU$_q$(1,1)
symmetry. As a result one can use the relevant SU$_q$(1,1) oscillator in
order to approximate the levels of the QESP which cannot be obtained using
that approach. In the specific case under discussion
the extra levels obtained in this way are the levels of the same parity lying
above the first $n+1$ known levels, as well as all the levels of the
opposite parity. The fact that the modified P\"oschl-Teller potential
with SU(1,1) symmetry cannot be connected to these QESPs, while its
deformed version with SU$_q$(1,1) symmetry can, is an example where the
use of q-deformation allows the solution of an otherwise intractable problem.

In ref. [12] the deformed anharmonic oscillator with SU$_q$(1,1) symmetry
has been proved appropriate for the description of vibrational spectra of
diatomic molecules. An extension to vibrational spectra of highly symmetric
polyatomic molecules has been given in [26]. The construction of QESPs
related to realistic vibrational spectra is receiving attention.

One of the authors (DB) has been supported by  CEC under contract
ERBCHBGCT 9\-3\-0\-4\-67. Another author
(HM) acknowledges support by KFUPM and the warm hospitality during his stay
in NCSR ``De\-mo\-kri\-tos''.

\vfill\eject

\centerline{\bf References}

\def\ref{\par\smallskip\hangindent=2cm\hangafter=1}
\parindent=0pt

\ref [1] V. G. Drinfeld, in {\it Proceedings of the International Congress of
Mathematicians}, ed. A. M. Gleason (American Mathematical Society, Providence,
RI, 1986) p. 798.

\ref [2] M. Jimbo, Lett. Math. Phys. 11 (1986) 247.

\ref [3] L. C. Biedenharn, J. Phys. A 22 (1989) L873.

\ref [4] A. J. Macfarlane, J. Phys. A 22 (1989) 4581.

\ref [5] M. Arik and D. D. Coon, J. Math. Phys. 17 (1976) 524.

\ref [6] V. V. Kuryshkin, Annales de la Fondation Louis de Broglie 5 (1980)
111.

\ref [7] A. Jannussis, in {\it Proceedings of the 5th International Conference
on Hadronic Mechanics}, ed. H. C. Myung (Nova Science, Commack, NY, 1990).

\ref [8] D. Bonatsos and C. Daskaloyannis, Phys. Lett. B 307 (1993) 100.

\ref [9] C. Daskaloyannis, J. Phys. A 24 (1991) L789.

\ref [10] D. Bonatsos, C. Daskaloyannis and K. Kokkotas, J. Phys. A 24
(1991) L795.

\ref [11] D. Bonatsos, C. Daskaloyannis and K. Kokkotas, J. Math. Phys.
33 (1992) 2958.

\ref [12] D. Bonatsos, E. N. Argyres and P. P. Raychev, J. Phys. A 24 (1991)
L403.

\ref [13] D. Bonatsos, P. P. Raychev and A. Faessler, Chem. Phys. Lett.
178 (1991) 221.

\ref [14] A. V. Turbiner and A. G. Ushveridze, Phys. Lett. A 126 (1987) 181.

\ref [15] A. V. Turbiner, Commun. Math. Phys. 118 (1988) 467.

\ref [16] A. Turbiner, Phys. Lett. B 276 (1992) 95.

\ref [17] C. M. Bender and A. Turbiner, Phys. Lett. A 173 (1993) 442.

\ref [18] A. Gonz\'alez-L\'opez, N. Kamran and P. J. Olver, Commun. Math.
Phys. 159 (1994) 503.

\ref [19] P. P. Kulish and E. V. Damaskinsky, J. Phys. A 23 (1990) L415.

\ref [20] H. Ui and N. Aizawa, Mod. Phys. Lett. A 5 (1990) 237.

\ref [21] P. Morse, Phys. Rev. 34 (1929) 57.

\ref [22] G. P\"oschl and E. Teller, Z. Phys. 83 (1933) 143.

\ref [23] I. H. McKenna and K. K. Wan, J. Math. Phys. 25 (1984) 1978.

\ref [24] Th. Ioannidou, Diploma Thesis, U. Thessaloniki (1993), unpublished.

\ref [25] C. P. Sun and W. Li, Commun. Theor. Phys. 19 (1993) 191.

\ref [26] D. Bonatsos and C. Daskaloyannis, Phys. Rev. A 48 (1993) 3611.

\vfill\eject

\centerline{\bf Table 1}

Exact energy levels of the quasi-exactly soluble potential of eq. (16)
with $n=9$, $r=0$, $a=1$, $b=21.275801$, compared to the approximate energy
levels of the WKB-EP potential of eq. (11) with $A=0.2703882$, $E_0'=
7126.0336$,
$\tau= 0.00545864$, $N=399$, given by eq. (9).

$$\vbox{\halign{\hfil #\hfil &&\quad \hfil #\hfil \cr
 $n$  & $E_n$ (approx) & $E_n$ (exact)  \cr
      &              &              \cr
0  & 20.3991 & 20.4153 \cr
2  & 102.672 & 102.682 \cr
4  & 186.131 & 186.138 \cr
6  & 270.735 & 270.749 \cr
8  & 356.444 & 356.485 \cr
10 & 443.217 & 443.317 \cr
12 & 531.013 & 531.218 \cr
14 & 619.791 & 620.165 \cr
16 & 709.507 & 710.133 \cr
18 & 800.118 & 801.101 \cr}}$$

\vfill\eject\bye